\documentclass[12pt, twoside]{article}

\usepackage{amsfonts}
\usepackage{amssymb}
\usepackage{graphics}
\usepackage[all]{xy}
\usepackage{amsmath, amsthm, latexsym}

\def \c{\mathbb{C}}
\def \z{\mathbb{Z}}
\def \r{\mathbb{R}}
\def \n{\mathbb{N}}
\def \p{\mathbb{P}}

\def \U{\mathcal{U}}

\def \A{\c[x_1^{\pm1}, \ldots,x_n^{\pm1}]}

\def \GL{\textup{GL}}
\def \SP{\textup{SP}}

\def \.{\cdot}

\def \dim{\textup{dim}}
\def \rank{\textup{rank}}

\def \det{\textup{det}}

\def \In{\textup{in}}

\def \Spec{\textup{Spec}}

\theoremstyle{plain}
\newtheorem{Th}{Theorem}[section]
\newtheorem{Lem}[Th]{Lemma}
\newtheorem{Prop}[Th]{Proposition}
\newtheorem{Cor}[Th]{Corollary}

\theoremstyle{definition}

\newtheorem{Def}[Th]{Definition}
\newtheorem{Rem}[Th]{Remark}

\begin{document}
\title{SAGBI bases and Degeneration of Spherical Varieties to Toric Varieties}
\author{Kiumars Kaveh\\Department of Mathematics\\University of British Columbia}
\maketitle

{\footnotesize {\bf Abstract.} Let $X \subset \p(V)$ be a
projective spherical $G$-variety, where $V$ is a finite
dimensional $G$-module and $G = \SP(2n, \c)$. In this paper, we
show that $X$ can be deformed, by a flat deformation, to the toric
variety corresponding to a convex polytope $\Delta(X)$. The
polytope $\Delta(X)$ is the polytope fibred over the moment
polytope of $X$ with the Gelfand-Cetlin polytopes as fibres. We
prove this by showing that if $X$ is a horospherical variety,
e.g. flag varieties and Grassmanians,
the homogeneous coordinate ring of $X$ can be
embedded in a Laurent polynomial algebra and has a SAGBI basis
with respect to a natural term order. Moreover, we show that the
semi-group of initial terms, after a linear change of variables,
is the semi-group of integral points in the cone over the polytope
$\Delta(X)$. The results of this paper are true for other
classical groups, provided that a result of A. Okounkov on the
representation theory of $\SP(2n,\c)$
is shown to hold for other classical groups.}\\


\tableofcontents

\section{Introduction} \label{sec-intro}
Let $X \subset \p(V)$ be a (normal) projective $G$-variety, where
$G$ is a classical group and $V$ is a finite dimensional
$G$-module. Suppose $X$ is {\it spherical}, that is a Borel
subgroup has a dense orbit. Generalizing the case of toric
varieties, one can associate an integral convex polytope
$\Delta(X)$ to $X$ such that the Hilbert polynomial $h(t)$ of $X$
is the Ehrhardt polynomial of $\Delta(X)$, i.e. $h(t) = $ number
of integral points in $t\Delta(X)$. The polytope $\Delta(X)$ is
the polytope fibred over the moment polytope of $X$ with the
Gelfand-Cetlin polytopes as fibres.
This polytope was
defined by A. Okounkov in ~\cite{Okounkov2},
based on the results of M. Brion.
We call this polytope the {\it Newton polytope} of $X$.

In this paper, for $G=\SP(2n,\c)$, we show that $X$ can be
deformed (degenerated), by a flat deformation, to the toric
variety corresponding to the polytope $\Delta(X)$ (Corollary
\ref{cor-main}). This is the consequence of the main result of the
paper, i.e. the homogeneous coordinate ring of a {\it
horospherical} variety has a SAGBI basis (Theorem \ref{thm-main}).
A spherical variety is horospherical if the stabilizer of a point
in the dense $G$-orbit contains a maximal unipotent subgroup. Flag
varieties and Grassmanians are examples of horospherical
varieties. It is known that any spherical variety can be deformed,
by a flat deformation, to a horospherical variety such that the
moment polytopes of the two varieties are the same (see
\cite{Popov}, \cite[$\S2.2$]{A-B}, \cite[Satz 2.3]{Knop}).

More precisely, we prove that if $X \subset \p(V)$ is a projective
horospherical $G$-variety where $G = \SP(2n,\c)$, the homogeneous
coordinate ring $R$ of $X$ can be embedded in a Laurent polynomial
algebra and has a SAGBI basis with respect to a natural term order
\footnote{SAGBI stands for {\it Subalgebra Analogue of Gr\"{o}bner
Basis for Ideals}.}. Moreover, we show that the semi-group of
initial terms is the semi-group of integral points in the cone
over the polytope $\Delta(X)$. A finite collection $f_1, \ldots,
f_r$ of elements of $R$ is a SAGBI basis, with respect to a term
order, if the semi-group of initial terms is generated by the
initial terms of the $f_i$ and moreover, every element of $R$ can
be represented as a polynomial in the $f_i$, in a finite number of
steps, by means of a simple classical algorithm called the {\it
subduction algorithm}.

Degenerations of flag and Schubert varieties to toric varieties have
been studied by Gonciulea and Lakshmibai in \cite{G-L}
and by Caldero in \cite{Caldero}. Recently,
M. Kogan and E. Miller show the existence of a SAGBI basis for the
coordinate ring of the flag variety of $\GL(n, \c)$. More precisely,
they prove that for any dominant weight $\lambda$ in the
interior of the Weyl chamber, the homogenous coordinate ring of the flag
variety $\GL(n)/B$ embedded in $\p(V_\lambda)$ has a SAGBI basis and
$\GL(n)/B$ can be degenerated to the toric variety corresponding to the
Gelfand-Cetlin polytope of $\lambda$ (see \cite{K-M}). Main results
of the present paper (Theorem \ref{thm-main} and Corollary \ref{cor-main})
, in particular, imply the similar result for
the flag varieties $G/P$ of $G=\SP(2n,\c)$.

A key step in our proof is a result of A. Okounkov on the
representation theory of $\SP(2n,\c)$. Let
$V_\lambda$ denote the irreducible $G$-module with highest weight
$\lambda$, where $G=\SP(2n,\c)$. It is well-known that one can
view $V_\lambda$ as a subspace of $\c[G]$ and, after restriction
to $U$, as a subspace of $\c[U]$, where $U$ is the standard
maximal unipotent subgroup of $G$. In \cite{Okounkov1}, Okounkov
proves that, with respect to a natural term order on $\c[U]$, the
set of highest terms of elements of $V_\lambda$ can be identified
with the Gelfand-Cetlin polytope $\Delta_\lambda$
(Theorem \ref{thm-Okounkov}). As Okounkov
informed the author, using similar methods used for $\SP(2n,\c)$,
one can prove his result for other classical groups. But so far he
has not published the proofs for other classical groups. The
results of the present paper as well as their proofs go verbatim
for other classical groups, provided that Okounkov's result
is shown to hold for them.

In Section \ref{sec-SAGBI}, we discuss SAGBI bases. Section
\ref{sec-homog-coordinate-ring} deals with some facts about
homogeneous coordinate ring of spherical varieties. We give a
description of the homogeneous coordinate ring of a horospherical
variety. In Section \ref{sec-Newton-polytope}, we define the
Gelfand-Cetlin polytopes and the polytope $\Delta(X)$. Section
\ref{sec-Okounkov} discusses the result of A. Okounkov on the
initial terms of elements of an irreducible $G$-module and
Gelfand-Cetlin polytopes, for $G = \SP(2n, \c)$. Finally, in
Section \ref{sec-main-theorem} we state and prove our main
results.\\

\noindent{\bf Acknowledgment:} The author would like to thank I.
Arzhantsev, J. Chipalkatti, A.G. Khovanskii, A. Okounkov and Z.
Reichstein for stimulating discussions. Also I would like to thank
I. Arzhantsev and Z. Reichstein and for reading the first version
and giving helpful comments.

\section{SAGBI bases}
\label{sec-SAGBI}
In this section we define the notion of a SAGBI basis for a
subalgebra of the Laurent polynomials. SAGBI bases play an important
role when one deals with subalgebras of the polynomial or Laurent
polynomial algebras. Their theory is more complicated than the theory of
Gr\"{o}bner bases. In particular, not every subalgebra has
a SAGBI basis with respect to a given term order. It is an
unsolved problem to determine, for a given term order, which subalgebras
have a SAGBI basis.

Let $\A$ denote the algebra
of Laurent polynomials in $n$ variables. Let $\prec$ be a term
order on $\z^n$, that is a total order compatible with addition.
An important example is the lexicographic order. The initial term,
with respect to $\prec$, of a polynomial $f$ is denoted by
$\In(f)$. If $R$ is a subalgebra of $\A$, we denote by $\In(R)$
the semi-group of initial terms in $R$, i.e. $\{ \In(f) \mid 0
\neq f \in R \}$.

First consider the case where $R$ is a subalgebra of
$\c[x_1,\ldots, x_n]$. In this case, one usually assumes that
$\prec$ satisfies the extra condition:
$$ {\bf a} \succ (0,\ldots,0),~~ \forall{\bf a}~~ 0\neq {\bf a} \in \n^n.$$

\begin{Def}
Let $R$ be a subalgebra of $\c[x_1,\ldots,x_n]$. A finite
collection of polynomials $\{ f_1, \ldots, f_r \} \subset R$ is a
SAGBI basis for $R$, if $\{\In(f_1), \ldots, \In(f_r)\}$ generates
the semi-group $\In(R)$.
\end{Def}

When $R$ has a SAGBI basis, one has a simple classical algorithm,
due to Kapur-Madlener and Robbiano-Sweedler, to express elements of $R$
in terms of the $f_i$ as follows: Write $\In(f) = d_1\In(f_1) +
\cdots + d_r\In(f_r)$ for some $d_1, \ldots, d_r \in \n$. Dividing
the leading coefficient of $f$ by the leading coefficient of
${f_1}^{d_1} \cdots {f_r}^{d_r}$, we obtain a $c$ such that the
leading term of $f$ is the same as the leading term of $c{f_1}^{d_1}
\cdots {f_r}^{d_r}$. Set $g = f - c{f_1}^{d_1} \cdots
{f_r}^{d_r}$. If $g = 0$, we are done; otherwise we replace $f$ by
$g$ and proceed inductively. Since $g$ has a smaller leading
exponent than $f$, and $\n^n$ is well-ordered with respect to
$\prec$, this process will terminate, resulting an expression for
$f$ as a polynomial in the $f_i$. This is referred to as {\it
subduction algorithm}. See \cite{Rob} for a detailed discussion of
SAGBI bases for subalgebras of $\c[x_1,\ldots,x_n]$.

In general when $R$ is a subalgebra of $\A$, since $\z^n$ is not
well-ordered there is no guarantee that this algorithm terminates.
Following \cite[p. 2]{Zinovy}, we define the SAGBI basis as
follows:
\begin{Def} \label{def-SAGBI}
Let $R$ be a subalgebra of $\A$. A finite collection of
polynomials $\{f_1, \ldots, f_r\}$ is a SAGBI basis for $R$ if:
\begin{itemize}
\item[(a)] The $\In(f_i)$ generate $\In(R)$ as a semi-group; and
\item[(b)] the subduction algorithm described above terminates for
every $f \in R$, no matter what choices are made for $d_1, \ldots,
d_r$ in the course of the algorithm.
\end{itemize}
The algebra $R$ is said to have a SAGBI basis, if it has a SAGBI
basis for some choice of a term order.
\end{Def}
\section{Homogeneous coordinate ring of spherical and horospherical varieties}
\label{sec-homog-coordinate-ring} Let $V$ be a finite dimensional
$G$-module and $X \subset \p(V)$ a projective spherical
$G$-variety, i.e. $X$ is normal and a Borel subgroup $B \subset G$ has a
dense orbit in $X$. Let $R = \c[X]$ denote the homogeneous coordinate ring of
$X$. This algebra is graded by the degree of polynomials, $$R =
\bigoplus_{k=0}^\infty R_k.$$ We decompose the spaces $R_k$ into
irreducible $G$-modules,
$$R_k = \bigoplus_{\lambda} m_{\lambda,k}V_\lambda,$$ where
$V_\lambda$ is the irreducible $G$-module with the highest weight
$\lambda$ and $m_{k,\lambda}$ is its multiplicity. Since $X$ is
spherical its spectrum is multiplicity free, i.e. $m_{k,\lambda}
\in \{0,1\}$. Let $\Phi(X)$ denote the {\it moment polytope} of
$X$, i.e. the intersection of the image of the moment map with the
positive Weyl chamber for the choice of $B$. Also, denote by
$\Lambda$ the weight lattice of $G$. The following theorem due to
Brion (see \cite{Brion1} and \cite{Brion2}) determines which
weights $\lambda$ occur in the decomposition of $R_k$ with
multiplicity $1$:
\begin{Th}[Brion, $\S3$ \cite{Brion2}] \label{thm-Brion}
There is a sublattice $\Lambda'$ of $\Lambda$ such that $\Phi(X)
\subset \Lambda'_\r$, the vector space spanned by $\Lambda'$, and
we have:
$$ R_k = \bigoplus_{\lambda \in k\Phi(X) \cap \Lambda'} V_\lambda.$$
\end{Th}
The rank of the sublattice $\Lambda'$ is called the {\it rank} of
the spherical variety $X$.

\begin{Rem} \label{rem-moment-polytope}
It follows from the above theorem that one can recover the moment
polytope $\Phi(X)$ from the multiplicities of the irreducible
$G$-modules appearing in $R_k$. More precisely, we have
$$\Phi(X) = \text{ closure of }\bigcup_{k=0}^{\infty} \{ \frac{\mu}{k}
\mid V_\mu
\text{ appears in the decomposition of } R_k \}.$$
\end{Rem}

One can show that the ring multiplication in $R$ sends $V_\lambda
\times V_\mu$ to $V_{\lambda + \mu} \oplus \bigoplus_\nu V_\nu$,
where $\nu = \lambda + \mu - \xi$ and $\xi$ is some non-negative
combination of simple roots. When all the stabilizer subgroups of
the points of $X$ contain a maximal unipotent subgroup, from a
theorem of Popov (see \cite[Theorem 2.3 ]{Popov}) it follows that
the ring multiplication sends $V_\lambda \times V_\mu$ to
$V_{\lambda + \mu}$ and this map coincides with a Cartan
multiplication. \footnote{For definition of Cartan multiplication
see \cite[p. 429]{F-H}}

\begin{Def} A spherical $G$-variety $X$
such that the stabilizer of a point in the dense
$G$-orbit contains a maximal unipotent subgroup is called a {\it
horospherical} variety.
\end{Def}
It can be shown that if $X$ is horospherical, then all the
stabilizer subgroups contain a maximal unipotent subgroup.
Examples of horospherical varieties are toric varieties, flag
varieties and Grassmanians.

Now, assume $X$ is horospherical. Fix a point $x$ in the dense
$G$-orbit of $X$. Choose highest weight vectors $f_\lambda$ in
each simple submodule $V_\lambda$ of $R$ by the condition that
$f_\lambda(x) = 1$. Then the product of these highest weight
vectors is again such a vector, and for any two $\lambda$ and
$\mu$ appearing in the decomposition of $R$, one can uniquely
define Cartan multiplication. We can then give the following
description for the homogeneous coordinate ring of $X$:
\begin{Th} \label{thm-coordinate-horo-var}
We have the following isomorphism of graded algebras:
$$ R \cong \bigoplus_{k=0}^{\infty}\bigoplus_{\lambda \in k\Phi(X) \cap \Lambda'}
V_\lambda,$$ where the multiplication in the righthand side is
defined as follows: Let $R_d = \bigoplus_{\lambda}V_\lambda$ and
$R_e = \bigoplus_\mu V_\mu$ be the decomposition of two graded
pieces of $R$. Then the multiplication $R_d \times R_e \to
R_{d+e}$ is given by the Cartan multiplication $V_\lambda \times
V_\mu \to V_{\lambda+\mu}$, defined uniquely by the above choice
of the highest weight vectors $f_\lambda$ and $f_\mu$.
\end{Th}

\section{Newton polytope of a spherical variety}
\label{sec-Newton-polytope} Let $G$ be a classical group. In this
section, following \cite{Okounkov2}, we briefly explain the
definition of the Newton polytope of a spherical $G$-variety $X$.
We start by recalling Gelfand-Cetlin polytopes.

To each dominant weight $\lambda$ of $G$, there corresponds a
Gelfand-Cetlin (or briefly G-C) polytope $\Delta_\lambda$.
The convex polytope $\Delta_\lambda$ has the property that the
number of integral points in $\Delta_\lambda$ is equal to the
dimension of the irreducible $G$-module $V_\lambda$. The dimension
of the Gelfand-Cetlin polytope is equal to the complex dimension
of the maximal unipotent subgroup $U$ of $G$, i.e.
$\frac{1}{2}(\dim(G) - \rank(G))$.
We recall the definition of Gelfand-Cetlin polytopes for $\GL(n, \c)$
and $\SP(2n, \c)$. For the definition of G-C polytopes for the
orthogonal group see \cite{B-Z}.
\begin{Def}
[G-C polytope for $\GL(n, \c)$]Let $\lambda = (\lambda_1 \geq
\cdots \geq \lambda_n)$ be a decreasing sequence of integers
representing a dominant weight in $\GL(n, \c)$. The G-C polytope
$\Delta_\lambda$ is the set of all real numbers $x_1, x_2, \ldots,
x_{n-1}$, $y_1, \ldots, y_{n-2}$, $\ldots, z$, such that the
following inequalities hold:
\begin{displaymath}
\begin{array}{ccccccccccc}
\lambda_1 && \lambda_2 && \lambda_3 & \cdots & \lambda_{n-2}
&& \lambda_{n-1} && \lambda_n \\
& x_1 && x_2 && \cdots && x_{n-2} && x_{n-1} & \\
&& y_1 && y_2 & \cdots & y_{n-3} && y_{n-2} &&\\
&&& \cdots && \cdots && \cdots\\ &&&\\
&&&& \cdots && \cdots && \\
&&&&& z &&&&&\\
\end{array}
\end{displaymath}
where the notation $$\begin{array}{ccc} a&&b\\&c& \end{array}$$
means $a \geq c \geq b$.
\end{Def}
\begin{Def}[G-C polytope for $\SP(2n, \c)$]
\label{def-G-C-SP2n} Let $B$ be the Borel subgroup of upper
triangular matrices in $\SP(2n, \c)$ and the maximal torus of
$\SP(2n, \c)$ be $\{ (t_1,\ldots,t_n,{t_1}^{-1}, \ldots, {t_n}^{-1})
\mid t_i \in \c^* , \forall i=1,\ldots,n \}$. Every dominant
weight is then represented by a decreasing sequence of positive
integers $\lambda = (\lambda_1 \geq \cdots \geq \lambda_n \geq
0)$. The G-C polytope $\Delta_\lambda$ is the set of all real
numbers $x_1, \ldots, x_{n}, y_1,\ldots,y_{n-1}, \ldots, z, w$,
such that the following inequalities hold:
\begin{displaymath}
\begin{array}{cccccccc}
\lambda_1 && \lambda_2 && \ldots & \lambda_n && 0 \\
& x_1 && x_2 & \ldots && x_n &\\
&& y_1 && \ldots & y_{n-1} && 0\\
&&& \ldots &&& \ldots &\\
&&&& \ldots &&&\\
&&&&& z && 0\\
&&&&&& w &\\
\end{array}
\end{displaymath}
\end{Def}

If the components of the weight $\lambda$ are real, we still can
define the $\Delta_\lambda$ by the above inequalities. So we can
extend the definition of $\Delta_\lambda$ to all real $\lambda$.

\begin{Lem}\label{Lem-G-C-is-linear}
The assignment $\lambda \mapsto \Delta_\lambda$ is linear, i.e.
$\Delta_{c\lambda} = c\Delta_\lambda$ for any positive $c$ and
$\Delta_{\lambda+\mu} = \Delta_\lambda + \Delta_\mu$, where the
addition in the righthand side is the Minkowski sum of convex
polytopes.
\end{Lem}
\begin{proof}
The proof is immediate from the definition in each of the three
cases of classical groups.
\end{proof}
Now, let $X \subset \p(V)$ be a (smooth) projective spherical
$G$-variety and $\Phi(X)$ its moment polytope. As before, let
$\Lambda$ denote the weight lattice and $\Lambda_\r$ the real
vector space spanned by $\Lambda$.
\begin{Def}[Newton polytope of a spherical variety]
Define the set $\Delta(X) \subset \Lambda_\r \oplus \r^{\dim U} =
\r^{\dim B}$, by $$ \Delta(X) = \bigcup_{\lambda\in \Phi(X)}
(\lambda, \Delta_\lambda).$$
\end{Def}
From Lemma \ref{Lem-G-C-is-linear}, it follows that $\Delta(X)$ is
a convex polytope.
\begin{Rem}\label{rem-Hilbert-poly}
In \cite{Okounkov2}, as a corollary of a theorem of Brion, it is
shown that the polytope $\Delta(X)$ has the property: $$\dim R_k =
\#\{k\Delta(X) \cap \Lambda' \},$$ where $\Lambda'$ is the
sublattice of $\Lambda$ in Theorem \ref{thm-Brion}. This means
that the Hilbert polynomial of the variety $X$ coincides with the
Ehrhardt polynomial of the polytope $\Delta(X)$. Note that
 since the Hilbert polynomial of a toric variety corresponding
to a polytope $\Delta$
is the Ehrhardt polynomial of $\Delta$, and the Hilbert polynomial
is invariant under a flat deformation, the above fact agrees with
the main result of the paper, i.e. $X$ can be deformed to the
toric variety of the polytope $\Delta(X)$ (Corollary \ref{cor-main}).
\end{Rem}

\section{Initial terms of elements of an irreducible $G$-module
and Gelfand-Cetlin polytopes}
\label{sec-Okounkov}
Let $\lambda$ be a dominant weight and
$V_\lambda$ the corresponding irreducible $G$-module, where $G=\SP(2n,\c)$.
The purpose of this section is to explain the result of A. Okounkov in
\cite{Okounkov1}, regarding the initial terms of the elements of
$V_\lambda$. We will need it in the proof of our main theorem.

First, we explain how one can identify $V_\lambda$ with a subspace
of a polynomial algebra, that is, the coordinate ring of the standard
maximal unipotent subgroup. Let $T$ be the standard maximal torus
of diagonal matrices in $G$,
$B_+$ the Borel subgroup of upper triangular matrices,
and $U_+$ the maximal unipotent
subgroup of $B_+$. Denote by $B_-$ and $U_-$ the opposite
subgroups of $B_+$ and $U_+$ respectively. Fix a $B_-$-eigenvector
$\xi$ in $(V_\lambda)^*$. It is well-known that the mapping from
$V_\lambda$ to $\c[G]$, defined by
$$ v \mapsto f_v,$$
$$ f_v(g) = \xi(g^{-1}v),$$
maps the $G$-module $V_\lambda$ isomorphically to the subspace
\begin{equation} \label{equ-1}
\{ f \in \c[G] \mid f(gb)=(-\lambda)(b)f(g) , \forall b\in B_- \}
\end{equation}
where $-\lambda$ is regarded as a character of $B_-$. We identify
$V_\lambda$ with its image in $\c[G]$. Choose the highest weight
vector $v_\lambda \in V_\lambda$ such that $\xi(v_\lambda) = 1$.

Consider the Bruhat decomposition $$G = \bigcup_{w \in W}
B_+wB_-,$$ where $W$ is the Weyl group. We have $G/B_- =
\bigcup_{w\in W} B_+wB_-/B_-$ and, the big Bruhat cell $\U$ in $G/B_-$
is $B_+B_-$. Since $B_+ \cap B_- = T$ and $B_+ = U_+T$, the cell
$\U$ can be identified with $U_+$, via $u \mapsto uB_-$. Since
$\U$ is dense in $G/B_-$, every element of $V_\lambda \subset
\c[G]$ is uniquely determined by its restriction to $U_+$. So we
can consider $V_\lambda$ as a subspace of $\c[U_+]$. Note that
$U_+$ is isomorphic, as a variety, to the affine space of
dimension $\frac{1}{2}(\dim(G) - \rank(G))$. One has:
\begin{Prop} \label{prop-comm-diag}
The following diagram is commutative: $$ \xymatrix{V_\lambda
\times V_\mu \ar@{^{(}->}[d] \ar[r] & V_{\lambda+\mu} \ar@{^{(}->}[d]\\
\c[G] \times \c[G] \ar[r]\ar[d] & \c[G] \ar[d]\\ \c[U_+] \times
\c[U_+] \ar[r] & \c[U_+]}$$ where the map in the first row is the
Cartan multiplication, defined uniquely with the above choice of
$v_\lambda$ and $v_\mu$, and the maps in the second and third rows
are the usual product of functions.
\end{Prop}
\begin{proof}
From (\ref{equ-1}) it follows that each $f_v$ defines a function
on $G/U_-$ and hence each $V_\lambda$ can be identified with a
subspace of $\c[G/U_-]$. Now the commutativity of the top part of
the diagram follows from a theorem of Popov (\cite[Theorem 2.
3]{Popov}, see also the paragraph after Remark
\ref{rem-moment-polytope}).
The commutativity of the bottom part of the diagram is
trivial.
\end{proof}
In \cite{Okounkov1}, Okounkov interprets
the G-C polytopes as the set of highest terms of the elements of
the $V_\lambda$ regarded as polynomials in $\c[U_+]$.
Choose a basis $e_1, \ldots, e_{2n}$ of $\c^{2n}$ in which the
matrix of the symplectic form is
$$\left[\begin{matrix} &&&&&1 \\ &0&&&\ldots& \\ &&&1&& \\ &&-1&&& \\
&\ldots&&&0& \\ -1&&&&& \\
\end{matrix}\right].$$

Let $x_{ij}$ be the matrix elements in this basis. We use $x_{11},\ldots,
x_{nn}$ as coordinates in $T$ and use the dual coordinates
$$g^\lambda = x_{11}^{\lambda_1}\cdots x_{nn}^{\lambda_n}, \quad g\in T, \lambda\in \Lambda,$$
for weights. The weights $$\lambda_1 \geq \lambda_2 \geq \cdots
\geq \lambda_n \geq 0$$ are dominant for $B_+$.

We use $x_{ij}, i<j, i+j \leq 2n+1$, as coordinates in $U_+$, and
the big Bruhat cell $\U$. Consider the following lexicographic
ordering on $\c[U_+]$: $$\prod x_{ij}^{p_{ij}} \succ \prod
x_{ij}^{q_{ij}} $$ if $p_{1,2n} < q_{1,2n}$, or if $p_{1,2n} =
q_{1,2n}$ and $p_{1,2n-1} < q_{1,2n-1}$, and so on. Note that in
particular
\begin{equation} \label{equ-2}
x_{1,2n} \prec x_{1,2n-1} \prec \cdots \prec x_{12} \prec
x_{2,2n-1} \prec \cdots \prec x_{23} \prec \cdots \prec x_{n,
n+1},
\end{equation}
which is exactly the reverse of the
ordering of positive roots induced by the standard lexicographic
order in $\r^n$. For a dominant weight $\lambda$ and a monomial
$$\prod x_{ij}^{p_{ij}},$$ put
\begin{eqnarray} \label{equ-3}
\eta_i &=& \lambda_i - p_{1,2n-i+1}, \quad i=1,\ldots,n, \cr
\theta_i &=& \eta_{i+1} + p_{1,i+1}, \quad i=1,\ldots,n-1, \cr
\eta'_i &=& \theta_i - p_{2,2n-i}, \quad i=1,\ldots,n-1, \cr
\theta'_i &=& \eta'_{i+1} + p_{2,i+1}, \quad i=1,\ldots,n-2,
\end{eqnarray}
\begin{Th}[\cite{Okounkov1}, Theorem 2] \label{thm-Okounkov}
View $V_\lambda$ as a subspace of $\c[U_+]$. Then, with the above
grading on $\c[U_+]$, the monomial $$\prod x_{ij}^{p_{ij}}$$ is a
highest monomial of a polynomial in $V_\lambda$ if and only if the
numbers $\eta_1,\ldots, \eta_n, \theta_1, \ldots, \theta_{n-1}
,\eta'_1, \ldots,\eta'_{n-1}, \ldots$, belong to the G-C polytope
$\Delta_\lambda$.
\end{Th}

Let us denote the vector $(\eta, \theta, \eta', \theta', \ldots) \in
\r^{\dim U}$ by $(q_{ij}),~ i<j, i+j\leq 2n+1$. The change of variables
$p_{ij} \mapsto q_{ij}$ in (\ref{equ-3}),
can be written in the matrix form as:
\begin{equation} \label{equ-4}
(q_{ij}) = A(p_{ij}) + B\lambda,
\end{equation}
where $A$ is a constant upper triangular matrix with $0, 1$ and
$-1$ as entries and $1, -1$ on the diagonal, and $B$ is the matrix
of the linear transformation
$$\lambda=(\lambda_1,\ldots,\lambda_n) \mapsto (\lambda_1,
\lambda_2 \ldots, \lambda_n, \lambda_2, \lambda_3
\ldots,\lambda_n,\ldots,\lambda_n) \in \r^{\dim(U)}.$$Note that
$\det(A) = \pm 1$ and hence the inverse of $A$ also has integer
entries. From (\ref{equ-4}) we can write
$$(p_{ij})=A^{-1}((q_{ij})-B\lambda),$$ Now, Theorem
\ref{thm-Okounkov} can be stated as follows: the monomial
$$\prod x_{ij}^{p_{ij}}$$
is a highest term of an element of $V_\lambda$ if and only if
$(p_{ij}) \in A^{-1}(\Delta_\lambda-B\lambda)$.
\begin{Def} \label{def-Delta'_lambda}
We denote the polytope $A^{-1}(\Delta_\lambda - B\lambda)$ by
$\Delta'_\lambda$.
\end{Def}
One has $\Delta_\lambda = A\Delta'_\lambda +B\lambda$, and hence
the two polytopes can be transformed to each other by integral
translations and integral transformations. Thus $\Delta_\lambda$
and $\Delta'_\lambda$ are integrally equivalent. The following is
immediate from the definition:
\begin{Lem}
The map $\lambda \mapsto \Delta'_\lambda$ is linear, i.e.
$\Delta'_{c\lambda} = c\Delta'_\lambda$ for a positive $c$, and
$\Delta'_{\lambda+\mu} = \Delta'_\lambda + \Delta'_\mu$ where the
addition in the righthand side is the Minkowski sum.
\end{Lem}
\begin{Def} \label{def-Delta'}
For a spherical variety $X$, similar to the definition of
$\Delta(X)$, define $\Delta'(X) \subset \Lambda_\r \oplus \r^{\dim
U} = \r^{\dim B}$, by
$$ \Delta'(X) = \bigcup_{\lambda\in \Phi(X)}
(\lambda, \Delta'_\lambda).$$
\end{Def}
From the above lemma, $\Delta'(X)$ is a convex polytope.

\begin{Rem} \label{rem-trans}
The map $(\lambda, x) \mapsto (\lambda, A^{-1}(x- B\lambda))$, is
an integral transformation that maps $\Delta(X)$ to $\Delta'(X)$.
The inverse of this transformation is $(\lambda, x) \mapsto
(\lambda, Ax+B\lambda)$ which is also integral. So the polytopes
$\Delta'(X)$ and $\Delta(X)$ can be transformed to each other by
integral transformations and hence are integrally equivalent.
\end{Rem}
\section{Main Theorem} \label{sec-main-theorem}
In this section, we prove the main results of the paper.
\begin{Th} \label{thm-main}
Let $V$ be a finite dimensional $G$-module,
and $X \subset \p(V)$ a projective horospherical $G$-variety, where $G=\SP(2n,\c)$.
We have:
\begin{itemize}
\item[(i)] The homogeneous coordinate ring $R$ of $X$ can be embedded into the
Laurent polynomial algebra
$\c[x_1,\ldots,x_d,y_1^{\pm1},\ldots,y_r^{\pm1},t]$, where
$d=\frac{1}{2}(\dim(G) - \rank(G))$ and $r=\rank(X)$.
\item[(ii)] $R$ has a SAGBI basis with
respect to a natural term order. Moreover, the semi-group of
initial terms $S = \In(R) \subset \z^{d+r+1}$ coincides with the
semi-group of integral points in the cone over the polytope
$\Delta'(X)$ (see Definitions \ref{def-Delta'_lambda} and
\ref{def-Delta'}), i.e.
$$ S = \z^{d+r+1} \cap \bigcup_{k=0}^\infty
(k\Delta'(X), k).$$
\end{itemize}
\end{Th}
\begin{proof}
We identify $\c[U_+]$ with the polynomial algebra
$\c[x_1,\ldots,x_d]$ equipped with the term order $\prec$ in
Theorem \ref{thm-Okounkov}. For each $\lambda$, let $\phi_\lambda$
denote the embedding ${V_{\lambda}} \hookrightarrow
\c[x_1,\ldots,x_d]$. Let $\Lambda'$ be the sublattice of the
weight lattice in Theorem \ref{thm-Brion}. Let $C \cong (\c^*)^r$
be a torus whose lattice of characters is $\Lambda'$. Let $y_1,
\ldots, y_r$ be a choice coordinates in $C$, hence $\c[C] =
\c[{y_1}^{\pm1},\ldots,{y_r}^{\pm1}]$. For $\lambda
=(\lambda_1,\ldots,\lambda_r) \in \Lambda'$, and $y =
(y_1,\ldots,y_r) \in C$, define $y^\lambda =
{y_1}^{\lambda_1}{y_2}^{\lambda_2}\ldots {y_r}^{\lambda_r}$.
Having the algebra isomorphism in Theorem
\ref{thm-coordinate-horo-var} in mind, define the function
$$\Psi: R = \bigoplus_{k=0}^\infty\bigoplus_{\lambda \in k\Phi(X)\cap\Lambda'}
V_{\lambda} \to
\c[x_1,\ldots,x_d,{y_1}^{\pm1},\ldots,{y_r}^{\pm1}, t],$$ by
$$\Psi(f) = t^k y^\lambda\phi_\lambda(f), \quad \forall f \in
V_{\lambda}, \lambda \in k\Phi(X)\cap\Lambda'$$ where $t$ is an
extra free variable. Then we have
\begin{Lem}
$\Psi$ is an injective homomorphism of algebras.
\end{Lem}
\begin{proof}
Since the $\phi_\lambda$ are additive homomorphisms, it
follows that $\Psi$ is also additive. The multiplicativity of
$\Psi$ follows from Proposition \ref{prop-comm-diag}.
$\Psi$ is 1-1, because the $\phi_\lambda$ are 1-1.
\end{proof}

Now, $R$ can be thought of as a subalgebra of $\c[x_1,\ldots,x_d,
{y_1}^{\pm1},\ldots,{y_r}^{\pm1}, t]$. Extend the term order
$\prec$ to $\c[x_1,\ldots,x_d,{y_1}^{\pm1},\ldots,{y_r}^{\pm1},t]$
by lexicographic order such that $t \succ y_r \succ \cdots \succ
y_1 \succ x_i, ~i=1,\ldots, d$. Let $S = \In(R) \subset
\z^{d+r+1}$. From Theorem \ref{thm-Okounkov}, we have
$$ S = \z^{d+r+1} \cap \bigcup_{k=0}^\infty\bigcup_{\lambda
\in k\Phi(X)\cap\Lambda'} (\Delta'_{\lambda}, \lambda, k),$$ i.e.
$S$ is the semi-group of integral points in the cone over the
polytope $\Delta'(X)$. This cone is a (strictly) convex rational
polyhedral cone and hence $S$ is finitely generated (Gordon's
lemma). Also, from the definition of $\prec$ and $S$, there are
only finitely many points in $S$ which are smaller than a given
point in $S$. This means that the subduction algorithm terminates
after a finite number of steps. Thus $R$ has a SAGBI basis and the
proof of the theorem is finished.
\end{proof}

Suppose $R$ is an arbitrary subalgebra of a Laurent polynomial
algebra. It is standard that the polynomials in $R$ can be
continuously deformed to their initial terms. More precisely, one
can show that there is a flat family of algebras $\pi: \mathcal{R}
\to \c$, such hat $\pi^{-1}(t) \cong R, \forall t \neq 0$ and
$\pi^{-1}(0) = \c[\In(R)]$, the semi-group algebra of $\In(R)$
(see \cite[Theorem 15.17]{Eisenbud}). If the semi-group $\In(R)$
is finitely generated then $\c[\In(R)]$ is the coordinate ring of
an affine (possibly non-normal) toric variety. Geometrically
speaking, this means that $\Spec(R)$ can be deformed, by a flat
deformation, to this affine toric variety.

\begin{Cor} \label{cor-horo}
Let $G=\SP(2n,\c)$.
Any projective horospherical $G$-variety $X \subset \p(V)$ can be
deformed, by a flat deformation, to the toric variety corresponding to the
polytope $\Delta(X)$. That is, there exists a flat family of
varieties $\pi: \mathcal{X} \to \c$, such that $\pi^{-1}(t) \cong X,
\forall t\neq 0$ and $\pi^{-1}(0)$ is the toric variety of the
polytope $\Delta(X)$.
\end{Cor}
\begin{proof}
Let $R$ be the homogeneous
coordinate ring of $X$. From \cite[Theorem 15.17, p.
343]{Eisenbud}, we know that $\Spec(R)$ can be deformed, by a flat
deformation, to the affine toric variety whose coordinate ring is
the semi-group algebra $\c[S]$. Since $\Delta'(X)$ and
$\Delta(X)$ can be transformed to each other by
integral transformations (Remark
\ref{rem-trans}), the semi-group $S$ is isomorphic to $S_0$, the
semi-group of integral points in the cone over $\Delta(X)$. So
$\Spec(R)$ can be deformed to the toric variety $\Spec(\c[S_0])$.
It is well-known that the projectivization of this affine toric
variety is the toric variety corresponding to the polytope
$\Delta(X)$ (see \cite{Sturmfels}, p. 36). This finishes the
proof of the corollary.
\end{proof}

Now, let $X \subset \p(V)$ be a projective spherical $G$-variety.
By a general result of Popov applied to the spherical varieties,
one can deform $X$, by a flat deformation, to a horospherical
variety $X_0$. More precisely:

\begin{Th} [see \cite{Popov}; \cite{A-B} $\S$2.2; \cite{Knop}
Satz 2.3]
Let $G$ be a reductive group and $Y$ an affine spherical $G$-variety.
There exists a flat family of affine $G$-varieties $\pi :\mathcal{Y} \to \c$
such that:
\begin{enumerate}
\item the $Y_t = \pi^{-1}(t)$ are isomorphic to $Y$ as $G$-varieties
for $t \neq 0$.
\item $Y_0 = \pi^{-1}(0)$ is horospherical.
\item $\c[Y]$ and $\c[Y_0]$ are isomorphic as graded $G$-modules,
in particular the multiplicities of the
irreducible representations $V_\lambda$ appearing in the graded pieces
$\c[Y]_d$ and $\c[Y_0]_d$ are the same, for any $d \geq 0$
\end{enumerate}
\end{Th}

If $X \subset \p(V)$ is a projective spherical variety, let $Y$ in the
above theorem be the cone over $X$ in $V$. We obtain that $X$ can be
degenerated to a projective horospherical variety $X_0$ where
$X_0$ is the projectivization of $Y_0$ in the theorem.
Since the multiplicities of the irreducible $G$-modules
apparing in the homogenuous
coordinate rings of $X$ and $X_0$ are the same we see that the
moment polytopes of $X$ and $X_0$ are the same
(see Remark \ref{rem-moment-polytope}). It is then
immediate from the definition that $\Delta(X) = \Delta(X_0)$.

\begin{Cor} \label{cor-main}
Let $G = \SP(2n, \c)$.
Any projective spherical $G$-variety
$X \subset \p(V)$ can be deformed,
by a flat deformation, to the toric variety corresponding to the
polytope $\Delta(X)$. That is, there exists a flat family of
varieties $\pi: \mathcal{X} \to \c$, such that $\pi^{-1}(t) \cong X,
\forall t\neq 0$ and $\pi^{-1}(0)$ is the toric variety of the
polytope $\Delta(X)$.
\end{Cor}
\begin{proof}
By the above comment $X$ can be deformed to a horospherical
variety $X_0$ and $\Delta(X) = \Delta(X_0)$.
The corollary now follows from Corollary \ref{cor-horo}.
\end{proof}

\noindent University of British Columbia, Vancouver, B.C. \\{\it
Email address:} {\sf kaveh@math.ubc.ca}

\end{document}